\documentclass{article}
\usepackage{graphicx}
\usepackage{psfrag}

\newtheorem{thm}{Theorem}[section]
\newtheorem{lem}{Lemma}[section]
\newtheorem{prop}{Proposition}[section]
\newtheorem{fact}{Fact}[section]
\newtheorem{qtn}{Question}[section]

\newenvironment{proof}{\medskip \noindent
{\bf Proof.}}{\hfill \rule{.5em}{1em} \\}
\newenvironment{pfmain}{\medskip \noindent
{\bf Proof of Theorem \ref{thm;main}.}}{\hfill \rule{.5em}{1em} \\}
\newenvironment{pfenough}{\medskip \noindent
{\bf Proof of Theorem \ref{thm;enough}.}}{\hfill \rule{.5em}{1em} \\}

\newcommand{\mb}{\mathbf}

\title{A family of pseudo metrics on $B^3$ and its application}  
\author{Young Deuk Kim \\ Department of Mathematics\\ State University of 
New York at Stony Brook\\ Stony Brook, NY 11794-3651, USA 
\\(ydkim@math.sunysb.edu)}
\date{\today}

\begin{document}
\maketitle

\begin{abstract}
Let $B^3$ be the closed unit ball in $\mb{R}^3$ and $S^2$ its boundary.
We define a family of pseudo metrics on $B^3$. As an application, we prove that
for any countable-to-one function $f:S^2\to [0,a]$, the set 
$$\mathcal{NM}^n_f=\{x\in S^2\mid \mbox{there exists } 
y\in S^2\ \mbox{such that } f(x)-f(y)>nd_E(x,y)\}$$ is uncountable 
for all $n\in\mb{N}$, where $d_E$ is the Euclidean metric on $\mb{R}^3$.  
\end{abstract}

\hspace{.5cm}
2000 Mathematics Subject Classification ; 57N05, 57M40

\section{The family of pseudo metrics}

In this section we construct the family of pseudo metrics on the closed unit
ball $B^3\subset\mb{R}^3$. As usual, a nonnegative function $d:B^3\times 
B^3\to\mb{R}$ is called a {\em pseudo metric} if 
\begin{enumerate}
\item
$d(x,x)=0$ for all $x\in B^3$
\item
$d(x,y)=d(y,x)$ for all $x,y\in B^3$
\item 
$d(x,y)+d(y,z)\geq d(x,z)$ for all $x,y,z\in B^3$.
\end{enumerate}

\indent
Let $S^2_r\subset\mb{R}^3$ be the 2-sphere with center $O=(0,0,0)$ and 
radius $0<r\leq 1$. We write $d_E$ to denote the Euclidean metric on 
$\mb{R}^3$. 
A metric $d$ on the set $S^2_r$ is called {\em locally Euclidean} if for all 
$P\in S^2_r$, there exists $t>0$ such that 
$$d(Q,R)=d_E(Q,R)\quad\mbox{for all }Q,R\in
B_t(P)=\{S\in S^2_r\mid d(P,S)<t\}.$$

\begin{figure}[ht]
\begin{center}
\psfrag{O}{$O$}
\psfrag{-P}{$-P$}
\psfrag{P}{$P$}
\psfrag{a}{$\alpha$}
\psfrag{b}{$2\alpha$}
\includegraphics[width=1.7in,height=1.8in]{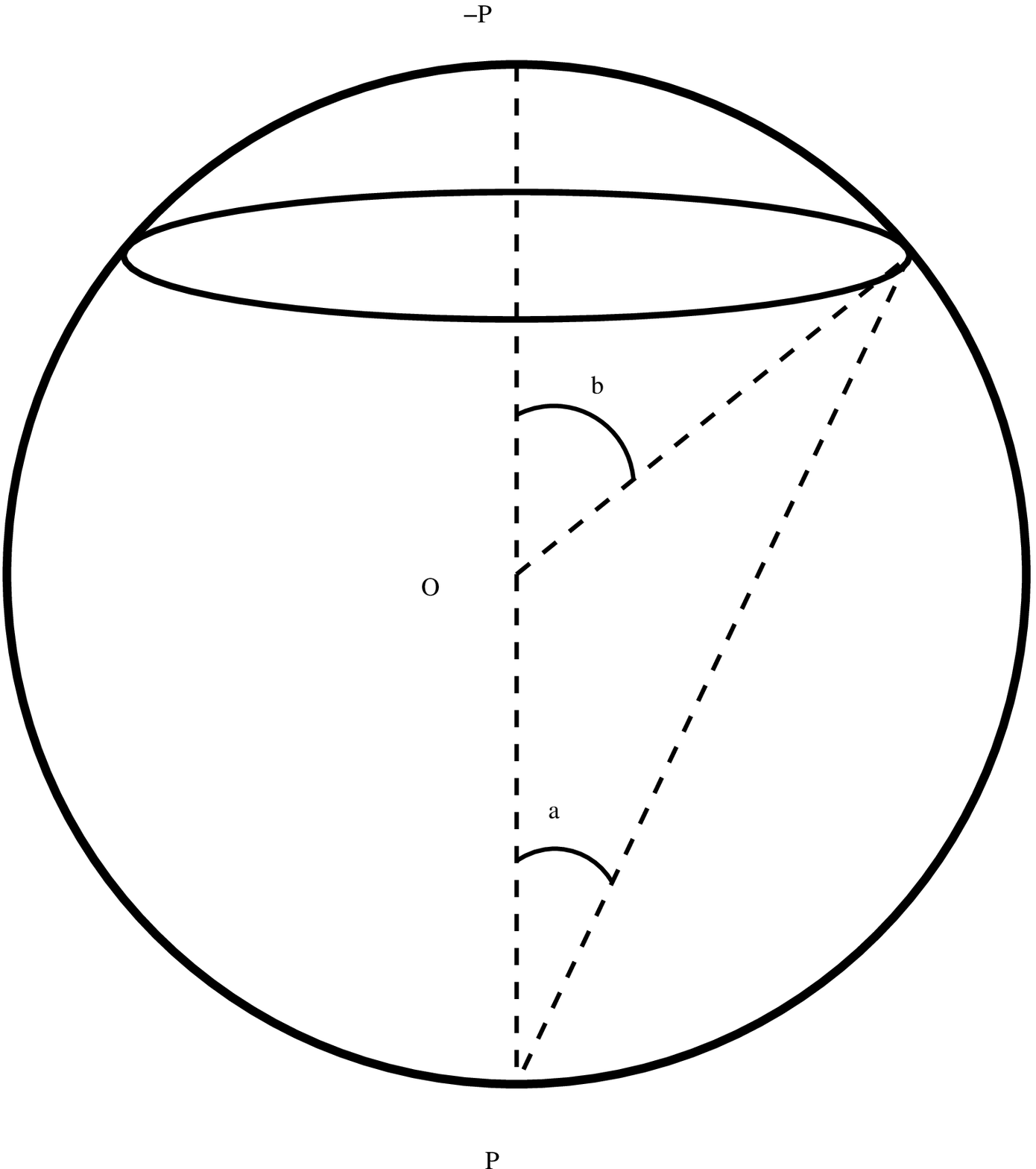}
\caption{$S^2_r$}
\label{fig;ds}
\end{center}
\end{figure}

\noindent
Suppose that $0<s\leq 1$. Let $-P$ denote the antipodal point of $P\in S^2_r$. 
Let 
$$\alpha=\sin^{-1}\left({\sqrt{2-s^2}-s}\over 2\right),\quad\mbox{where }
0\leq\alpha<\pi/4.$$ 
Notice that $\sin\alpha$ is a decreasing function of $s$ and hence so is 
$\alpha$. For all $P,Q\in S^2_r$, let (see Figure \ref{fig;ds}) 
\begin{eqnarray*}
d^s(P,Q)=\left\{
\begin{array}{lll}
d_E(P,Q) &\mbox{ if }& \angle POQ\leq \pi-2\alpha\\ 
2rs+d_E(-P,Q) &\mbox{ if }& \angle POQ>\pi-2\alpha,  
\end{array}
\right.
\end{eqnarray*}

\noindent
where $\alpha$ is defined as above.
Notice that if $s=1$ then $d^1=d_E$, and $d^s(P,-P)=2rs$ for all $P\in S^2_r$.
In \cite{KIM} the author proved 
\begin{thm} 
For all $0<s\leq 1$, $d^s$ is a locally Euclidean metric on $S^2_r$ which is
invariant under any Euclidean isometry.
\end{thm}

\noindent
Note that (see \cite{KIM} for a proof)
\begin{equation}\label{eqn;01}
sd_E(P,Q)\leq d^s(P,Q)\leq d_E(P,Q)\quad\mbox{for all }P,Q\in S^2_r
\end{equation}
and if $d^s(P,Q)\neq d_E(P,Q)$ then 
\begin{equation}\label{eqn;02}
 d^s(-P,Q)=d_E(-P,Q).
\end{equation}

\begin{figure}[ht]
\begin{center}
\psfrag{P}{$P$}
\psfrag{r}{$r$}
\psfrag{S1}{$S^2_r$}
\psfrag{S2}{$S_\lambda$}
\psfrag{B}{$B^3$}
\psfrag{E}{$\epsilon_\lambda$}
\includegraphics[width=3in,height=2in]{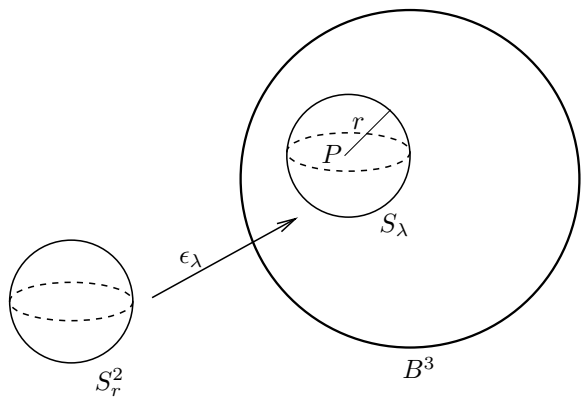}
\caption{$\lambda=(P,r)\in\Delta$}
\label{fig;slambda}
\end{center}
\end{figure}

\indent
We are going to use the metric spaces $(S^2_r,d^s)$ with $0<r,s\leq 1$ 
to construct a family of pseudo metrics on $B^3$. Let (see Figure 
\ref{fig;slambda}) $$\Delta=\left\{(P,r)\mid P\in int(B^3)\ \mbox{and } 
0<r\leq1-d_E(O,P)\right\}.$$ 

\noindent
For each $\lambda=(P,r)\in\Delta$, choose and fix an isometric embedding 
$\epsilon_\lambda$ of $\left(S^2_r,d_E\right)$ into $(B^3,d_E)$ such that 
$$S_\lambda=\epsilon_\lambda\left(S^2_r\right)
=\left\{Q\in B^3\mid  d_E(P,Q)=r\right\}.$$
Suppose that $\Lambda\subset\Delta$ and $P,Q\in B^3$. Let
\begin{eqnarray*}
&\Pi^n=(B^3)^{n+1}\times {\Lambda}^n\\ 
&[X,\lambda,n]=(X_0,X_1,\cdots,X_n,\lambda_1,\cdots,\lambda_n)\in\Pi^n \\ 
&\Gamma_{P,Q}=\left\{ [X,\lambda,n]\in\Pi^n\mid X_0=P,\ X_n=Q,\  
\mbox{and }X_{i-1},X_i\in S_{\lambda_i}\ \mbox{for all }1\leq i\leq n\right\}.
\end{eqnarray*}

\indent
Let $\Omega=\{\Lambda\subset\Delta\mid\Gamma_{P,Q}\neq\emptyset\ 
\mbox{for all }P,Q\in B^3\}$. 
Suppose that $\Lambda\in\Omega$ and $s:\Lambda\to (0,1]$ is a function. 
Define a function $d^{\Lambda,s}:B^3\times B^3\rightarrow\mb{R}$ by
$$d^{\Lambda,s}(P,Q)=\inf_{[X,\lambda,n]\in \Gamma_{P,Q}} \sum^n_{i=1}
d^{s(i)}\left({\epsilon_i}^{-1}(X_{i-1}),
{\epsilon_i}^{-1}(X_i)\right),$$

\noindent
where $s(i)=s(\lambda_i)$ and $\epsilon_i=\epsilon_{\lambda_i}$.
Notice that $d^{s(i)}\left({\epsilon_i}^{-1}(X_{i-1}),{\epsilon_i}^{-1}(X_i)
\right)$ does not depend on $\epsilon_i$ because $d^s$ is 
invariant under Euclidean isometries for all $s$.

\indent
Suppose that $P,Q\in S_\lambda\subset B^3$ for some $\lambda\in
\Delta$. Throughout this paper, we write $d^{s(\lambda)}(P,Q)$ to denote
$d^{s(\lambda)}\left(\epsilon_\lambda^{-1}(P),\epsilon_\lambda^{-1}(Q)\right)$.
 Therefore we have  
$$d^{\Lambda,s}(P,Q)=\inf_{[X,\lambda,n]\in \Gamma_{P,Q}} \sum^n_{i=1}
d^{s(i)}(X_{i-1},X_i).$$ 

\noindent
We also write $d_E(P,Q)$ to denote $d_E\left(\epsilon_\lambda^{-1}(P),
\epsilon_\lambda^{-1}(Q)\right)$. Notice that there is no problem with this 
notation because $\epsilon_\lambda$ is an isometric embedding for all
$\lambda\in\Delta$.

\indent
It is straightforward that $d^{\Lambda,s}$ is a pseudo metric on 
$B^3$ for all $\Lambda\in\Omega$ and $s:\Lambda\to (0,1]$.

\section{Elementary properties of the pseudo metric} 

In this section we study some elementary properties of the pseudo metric.
In particular we are interested in the following question.

\begin{qtn}
When is $d^{\Lambda,s}$ a metric on $B^3$? 
If $d^{\Lambda,s}$ is a metric on $B^3$, when is $(B^3,d^{\Lambda,s})$
homeomorphic to $(B^3, d_E)$?
\end{qtn}

\indent
Suppose that $\Lambda\in\Omega$, $s:\Lambda\to (0,1]$ and $P,Q\in B^3$. 
The following subset of $\Gamma_{P,Q}$ will turn out to be useful.
\begin{eqnarray*}
&&\Gamma^{\mathcal A}_{P,Q}=\big\{ [X,\lambda,n]\in\Gamma_{P,Q}\mid 
\mbox{for all }1\leq i\leq n,\\  
&&\qquad X_{i-1}\ \mbox{and }X_i\ \mbox{are antipodal on }S_{\lambda_i}\  
\mbox{or }\  d^{s(i)}(X_{i-1},X_i)=d_E(X_{i-1},X_i)\big\} 
\end{eqnarray*}

\indent
Let $[X,\lambda,n]=(X_0,X_1,\cdots,X_n,\lambda_1,\cdots,\lambda_n)
\in\Gamma_{P,Q}$ and $\lambda_i=(P_i,r_i)$.
Suppose that $d^{s(i)}(X_{i-1},X_i)\neq d_E(X_{i-1},X_i)$ for some 
$1\leq i\leq n$.
Let $-X_{i-1}$ be the antipodal point of $X_{i-1}$ on $S_{\lambda_i}$.
Notice that $$\epsilon_i^{-1}(-X_{i-1})=-\epsilon_i^{-1}
(X_{i-1}),$$ where $-\epsilon_i^{-1}(X_{i-1})$ is the antipodal point of
$\epsilon_i^{-1}(X_{i-1})$ on $S^2_{r_i}$. From eq. (\ref{eqn;02}) 
we have $d^{s(i)}(-X_{i-1},X_i)=d_E(-X_{i-1},X_i)$. Therefore 
\begin{eqnarray*}
&&d^{s(i)}(X_{i-1},X_i) \\
&=&2r_is(\lambda_i)+d_E(-X_{i-1},X_i)\\
&=&d^{s(i)}(X_{i-1},-X_{i-1})+d^{s(i)}(-X_{i-1},X_i).
\end{eqnarray*}

\noindent
From this observation we have
\begin{equation}\label{eqn;calk}
\qquad d^{\Lambda,s}(P,Q)=\inf_{[X,\lambda,n]\in 
\Gamma^{\mathcal A}_{P,Q}}\sum^n_{i=1}d^{s(i)}(X_{i-1},X_i). 
\end{equation}

\begin{figure}[ht]
\begin{center}
\psfrag{P}{$P$}
\psfrag{Q}{$Q$}
\psfrag{1}{$X_1$}
\psfrag{2}{$X_2$}
\psfrag{n-1}{$X_{n-1}$}
\psfrag{S1}{$S_{\lambda_1}$}
\psfrag{S2}{$S_{\lambda_2}$}
\psfrag{Sn}{$S_{\lambda_n}$}
\includegraphics[width=2.7in,height=1.2in]{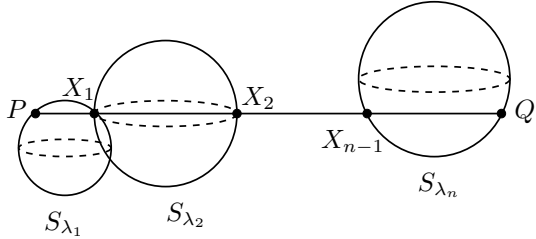}
\caption{piecewise dense}
\label{fig;pd}
\end{center}
\end{figure}

\indent
Suppose that $\Lambda\in\Omega$. $\Lambda$ is called {\em piecewise dense} if 
for all distinct two points $P,Q\in B^3$, there exists 
$[X,\lambda,n]\in\Gamma_{P,Q}$ (see Figure \ref{fig;pd})
such that $X_i$ is on the Euclidean segment $\overline{PQ}$ and 
$d_E(P,X_{i-1})<d_E(P,X_i)$ for all $1\leq i\leq n$.

\noindent
Note that, for these $X_i$'s, we have $\sum^n_{i=1}d_E(X_{i-1},X_i)=d_E(P,Q)$.
Let $$\Delta_t=\left\{(P,r)\in\Delta\mid r\leq t\right\}.$$
Notice that if $\Delta_t\subset\Lambda$ for some $t>0$ then 
$\Lambda$ is piecewise dense.

\indent
By the following theorem if $\Lambda$ is piecewise dense and 
$s$ is bounded below by some $s_0>0$, then $d^{\Lambda,s}$ is a metric on 
$B^3$ and the identity map from $(B^3,d_E)$ to $(B^3,d^{\Lambda,s})$ is a 
homeomorphism.
   
\begin{thm}\label{thm;ch}
Suppose that $\Lambda\in\Omega$, $s:\Lambda\rightarrow (0,1]$ and $P,Q\in B^3$.
\begin{enumerate}
\item
If $\Lambda$ is piecewise dense then $d^{\Lambda,s}(P,Q)\leq d_E(P,Q)$.
\item
If $s$ is bounded below by some $s_0>0$ then 
$d^{\Lambda,s}(P,Q)\geq s_0d_E(P,Q)$. 
\end{enumerate}
\end{thm}

\begin{proof}
Suppose that $\Lambda$ is piecewise dense and $P,Q\in B^3$.
We may assume that $P\neq Q$. Choose $[X,\lambda,n]\in\Gamma_{P,Q}$ such that
$X_i$ is on the Euclidean segment $\overline{PQ}$ and  $d_E(P,X_{i-1})<
d_E(P,X_i)$ for all $1\leq i\leq n$. From eq. (\ref{eqn;01}) we have 
$$d^{\Lambda,s}(P,Q)\leq \sum^n_{i=1}d^{s(i)}(X_{i-1},X_i)\leq
\sum^n_{i=1}d_E(X_{i-1},X_i)=d_E(P,Q).$$ 

\indent
Suppose that $s$ is bounded below by some $s_0>0$ and $P,Q\in B^3$.
From eq. (\ref{eqn;01}) we have 
\begin{eqnarray*}
d^{\Lambda,s}(P,Q)&=&\inf_{[X,\lambda,n]\in \Gamma_{P,Q}}\sum^n_{i=1}
d^{s(i)}(X_{i-1},X_i) \\
&\geq& \inf_{[X,\lambda,n]\in\Gamma_{P,Q}}\sum^n_{i=1}s(\lambda_i)
d_E(X_{i-1},X_i) \\
&\geq& \inf_{[X,\lambda,n]\in\Gamma_{P,Q}}\sum^n_{i=1}s_0d_E(X_{i-1},X_i) \\
&=&s_0\inf_{[X,\lambda,n]\in\Gamma_{P,Q}}\sum^n_{i=1}d_E(X_{i-1},X_i) \\
&\geq& s_0d_E(P,Q). 
\end{eqnarray*}
\end{proof}

\section{The application}

The application in this section is motivated by the famous 3-dimensional
Poincar\'{e} conjecture.

\vspace{.3cm}
\noindent
{\bf Poincar\'e Conjecture}
{\em If a compact connected 3-manifold is homotopic to $S^3$ then it is 
homeomorphic to $S^3$.}
\vspace{.3cm}

\noindent
See \cite{ANDERSON,HAMILTON,MILNOR2,PERELMAN1,PERELMAN2,PERELMAN3}
for the recent progress on the Poincar\'{e} conjecture by Perelman.
The work of Perelman is geometric and analytic. 
See \cite{BING,GR,HAKEN,MILNOR1,POENARU,ROURKE1,ROURKE2}  
for the topological approach to the Poincar\'{e} conjecture. 

\indent
Let $(M,d)$ be a compact connected 3-manifold with metric $d$ which is 
homotopic to $S^3$. Suppose that $h:(M,d)\to (S^3,d_E)$ be the homotopy 
equivalence, where $d_E$ is the Euclidean metric on $S^3\subset\mb{R}^4$.
Let $g:(S^3,d_E)\to (M,d)$ be any function such that $g\circ h(m)=m$ for all
$m\in M$. Note that $g$ is not necessarily continuous.

\begin{qtn}\label{qtn;01}
Suppose that $x,y\in S^3$. How can we compare $d(g(x),g(y))$ with $d_E(x,y)$?
\end{qtn}

\indent
One possible approach to the above question is as follows. 
For a fixed point $m_0\in M$, define a function $f:(S^3,d_E)\to\mb{R}$ 
by $f(x)=d(m_0,g(x))$. Notice that $d(g(x),g(y))\geq f(x)-f(y)$ for all 
$x,y\in S^3$ and $f(S^3)\subset [0,a]$ for some $a>0$ because $(M,d)$ is 
compact. Therefore the following question could be a subquestion of
Question \ref{qtn;01}.

\begin{qtn}
Suppose that $f:(S^3,d_E)\to [0,a]$ is a function and $x,y\in S^3$. 
What can we say about $f(x)-f(y)$ and $d_E(x,y)$?
\end{qtn}

\indent
There is no countable-to-one continuous function from $(S^3,d_E)$ to the 
closed interval $[0,a]$. And any countable-to-one function on $S^3$ induces
countable-to-one functions on subsets of $S^3$. 
As an application of the pseudo metric of previous sections, we will prove 
the following theorem in the next section.

\begin{thm}[Main Theorem]\label{thm;main}
Let $S^2$ be the boundary of $B^3\subset\mb{R}^3$.
Suppose that $a>0$ and $f:S^2\to [0,a]$ is a countable-to-one function.
Then $\mathcal{NM}^n_f$ is uncountable for all $n\in\mb{N}$, where 
$$\mathcal{NM}^n_f=\left\{x\in S^2\mid\mbox{there exists }y\in S^2\ 
\mbox{such that }f(x)-f(y)>nd_E(x,y)\, \right\}.$$
\end{thm}
 
\indent
Note that $\mathcal{NM}^n_f $ does not contain the set of 
discontinuities of $f$. For example, if $\mathit{C}$ be a countable dense 
subset of $S^2$ and 
\[ f(x)=\left\{ \begin{array}{rl}
             1& \mbox{ if }\; x\in \mathit{C} \\
             0&  \mbox{ otherwise}
           \end{array}\right.   \] 
then $\mathcal{NM}^n_f=\mathit{C}$ for all $n\in\mb{N}$ but $f$ is 
discontinuous at every point on $S^2$. 

\indent
Notice that Theorem \ref{thm;main} is not true for 
$S^1=\{(x,y)\mid x^2+y^2=1\}\subset\mb{R}^2$.
For example, if $g(x)=d_E(x,x_0)$ for some fixed $x_0\in S^1$ then
$g^{-1}(t)$ has at most two elements for all $t$ but $\mathcal{NM}^1_g=
\emptyset$ because 
$$g(x)-g(y)=d_E(x,x_0)-d_E(y,x_0)\leq d_E(x,y)\quad
\mbox{for all }x,y\in S^1.$$

\noindent
Notice also that countability of $f^{-1}(t)$ is crucial in Theorem 
\ref{thm;main}. For example, if $h(x)=d_E(x,x_0)$ for some fixed 
$x_0\in S^2$ then $h^{-1}(t)$ is uncountable for all $0<t<2$ and 
$\mathcal{NM}^1_h=\emptyset$.

\section{Proof of Theorem \ref{thm;main}}

In this section we prove Theorem \ref{thm;main}. We write $S^2$ to denote the 
boundary of $B^3\subset\mb{R}^3$ throughout this section. First, we will show 
that the following theorem implies Theorem \ref{thm;main}.

\begin{thm}\label{thm;enough}
Suppose that $f:S^2 \rightarrow [0,1]$ is a countable-to-one function.  
Then $\mathcal{NM}^2_f$ is uncountable.
\end{thm}

\indent 
Theorem \ref{thm;enough} implies the following theorem.

\begin{thm}\label{thm;a}
Suppose that $a>0$ and $f:S^2\rightarrow [0,a]$ is a countable-to-one function.
Then $\mathcal{NM}^2_f$ is uncountable.
\end{thm}

\begin{proof}
Choose $n\in \mb{N}$ such that $n>a$ and define a function $g$ by 
$g(x)={1\over n}f(x)$. Then 
$g:S^2\rightarrow \left[0,{a\over n}\right ]\subset[0,1]$ and  
\begin{eqnarray*}
\mathcal{NM}^2_f&=&\{x\in S^2\mid \exists y\in S^2\ \mbox{such that } 
f(x)-f(y)>2d_E(x,y)\} \\
&=&\{x\in S^2\mid \exists y\in S^2\ \mbox{such that } 
ng(x)-ng(y)>2d_E(x,y)\} \\
&\supset& \{x\in S^2\mid \exists y\in S^2\ \mbox{such that } 
g(x)-g(y)>2d_E(x,y)\} \\
&=&\mathcal{NM}^2_g. 
\end{eqnarray*}
Notice that $g^{-1}(t)$ is countable for all $t$ because 
$g^{-1}(t)=f^{-1}(nt)$. Therefore by Theorem \ref{thm;enough}, 
$\mathcal{NM}^2_g$ is uncountable. Hence $\mathcal{NM}^2_f$ is uncountable.  
\end{proof}

\indent
Due to Theorem \ref{thm;a}, we can prove Theorem \ref{thm;main}.

\begin{pfmain}
Define a function $h$ by $h(x)={1\over n}f(x)$. Notice that $\mathcal{NM}^n_f
\supset\mathcal{NM}^2_h$. Notice also that $h^{-1}(t)$ is countable for all 
$t$ because $h^{-1}(t)=f^{-1}(nt)$. Therefore by Theorem \ref{thm;a}, 
$\mathcal{NM}^2_h$ is uncountable. Hence $\mathcal{NM}^n_f$ is uncountable.  
\end{pfmain}

\begin{figure}[ht]
\begin{center}
\psfrag{x}{$x$}
\psfrag{S1}{$S_{({x\over 2},{1\over 2})}$}           
\psfrag{S2}{$S_{(P,{1\over 2})}$} 
\psfrag{S3}{$S_{(P,r)}$} 
\includegraphics[width=2.2in,height=2.2in]{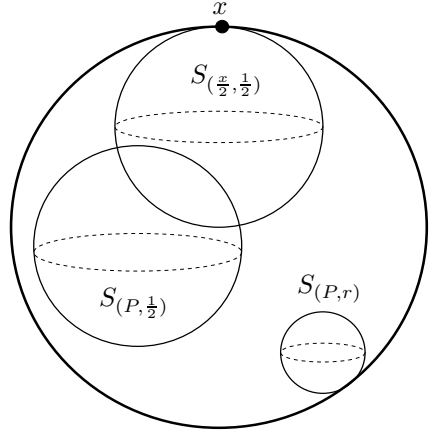}
\caption{$s(P,r)=1$ if $r\neq {1\over 2}$ or $S_{(P,r)}$ does not intersect 
$S^2$.}
\label{fig;application}
\end{center}
\end{figure}

\indent
In the remaining of this paper we will prove Theorem \ref{thm;enough}.
Suppose that $f:S^2\rightarrow [0,1]$ is a function. Notice that if 
$x\in S^2$ then ${x\over 2}\in int(B^3)$. From now on, 
we will use the following fixed $\Lambda\in\Omega$ and $s:\Lambda\to (0,1]$ 
to prove Theorem \ref{thm;enough} (see Figure \ref{fig;application}).
\begin{eqnarray*} 
&&\Lambda=\Delta_{1\over 2}=
\left\{(P,r)\in\Delta \mid r\leq {1\over 2} \right\}\\
 &&s\left({x\over 2},{1\over 2}\right)={{1+f(x)}\over 2}\ \ 
\mbox{for all } 
x\in S^2\quad \mbox{and}\quad s=1\  \mbox{otherwise.}
\end{eqnarray*}  

\noindent
Notice that $\Lambda$ is piecewise dense and the function $s$ is bounded below 
by ${1\over 2}$. By Theorem \ref{thm;ch}, $d^{\Lambda,s}$ is a metric
on $B^3$ and the identity map from $(B^3,d_E)$ to 
$(B^3,d^{\Lambda,s})$ is homeomorphism. 

\indent
The following two lemmas will turn out to be useful.

\begin{lem}\label{lem;half}
If $d^{\Lambda,s}(P,Q)<{1\over 2}$ then   
$d^{\Lambda,s}(P,Q)=d_E(P,Q)$ for all $P,Q\in B^3$. 
\end{lem}

\begin{proof}
Suppose that $P,Q\in B^3$ and $d^{\Lambda,s}(P,Q)<{1\over 2}$. Let 
$$\Gamma^E_{P,Q}=\left\{ [X,\lambda,n]\in\Gamma^{\mathcal A}_{P,Q}
\mid \mbox{for all }1\leq i\leq n,\ 
d^{s(i)}(X_{i-1},X_i)=d_E(X_{i-1},X_i)\right\}.$$
Suppose that $[X,\lambda,n]\in\Gamma^{\mathcal A}_{P,Q}\setminus
\Gamma^E_{P,Q}$. There exists $j$ such that 
$X_{j-1}$ and $X_j$ are antipodal on $S_{\lambda_j}$ with $\lambda_j=
\left({x\over 2},{1\over 2}\right)$ for some $x\in S^2$. Therefore
$$\sum^n_{i=1}d^{s(i)}(X_{i-1},X_i)\geq d^{s(j)}(X_{j-1}
,X_j)={{1+f(x)}\over 2}\geq {1\over 2}.$$ 

\indent
Recall that $d^{\Lambda,s}(P,Q)<{1\over 2}$. Therefore 
from eq. (\ref{eqn;calk}), we have
\begin{eqnarray*}
d^{\Lambda,s}(P,Q)&=&\inf_{[X,\lambda,n]\in \Gamma^{\mathcal A}_{P,Q}}
\sum^n_{i=1}d^{s(i)}(X_{i-1},X_i) \\
&=&\inf_{[X,\lambda,n]\in \Gamma^E_{P,Q}}
\sum^n_{i=1}d^{s(i)}(X_{i-1},X_i) \\
&=&\inf_{[X,\lambda,n]\in \Gamma^E_{P,Q}}
\sum^n_{i=1}d_E(X_{i-1},X_i) \\
&\geq& d_E(P,Q).
\end{eqnarray*}
Since $\Lambda$ is piecewise dense, we have 
$d^{\Lambda,s}(P,Q)\leq d_E(P,Q)$. Therefore
$$d^{\Lambda,s}(P,Q)=d_E(P,Q).$$   
\end{proof}

\begin{lem}\label{lem;haf}
If $d_E(P,Q)<{1\over 2}$ then $d^{\Lambda,s}(P,Q)=d_E(P,Q)$ for all 
$P,Q\in B^3$. 
\end{lem}

\begin{proof}
Suppose that $d_E(P,Q)<{1\over 2}$.
Since $\Lambda$ is piecewise dense, we have 
$d^{\Lambda,s}(P,Q)\leq d_E(P,Q)<{1\over 2}$.
Therefore from Lemma \ref{lem;half}, we have $d^{\Lambda,s}(P,Q)=d_E(P,Q)$.   
\end{proof}

\indent
Suppose that $O=(0,0,0)\in B^3$ and $x\in S^2$. 
Recall from eq. (\ref{eqn;calk}) that
\begin{equation}
d^{\Lambda,s}(O,x)=\inf_{[X,\lambda,n]\in\Gamma^{\mathcal A}_{O,x}}
\sum^n_{i=1}d^{s(i)}(X_{i-1},X_i). \label{eqn;kdis} 
\end{equation}

\noindent
The following subsets of $\Gamma^{\mathcal A}_{O,x}$ will turn out to be 
useful.
\begin{eqnarray*}
\Gamma^{\mathcal U}_{O,x}&=&\Big\{ [X,\lambda,n]\in\Gamma^{\mathcal A}_{O,x}
\mid \mbox{there exists unique } i\ \mbox{such that }
X_{i-1}\ \mbox{and }X_i\ \\
&&\qquad \mbox{are antipodal on }S_{\lambda_i}\
\mbox{with }\lambda_i=\left({y\over 2},{1\over 2}\right)\ 
\mbox{for some } y\in S^2 \Big\} \\ 
\Gamma^{\partial}_{O,x}&=&\Big\{ [X,\lambda,n]\in \Gamma^{\mathcal U}_{O,x}
\mid \mbox{there exists }y\in S^2\ \mbox{such that }\\ 
&&\qquad\qquad\qquad X_1=y\ \mbox{and }\lambda_1=\left({y\over 2},{1\over 2}
\right)\Big\}. 
\end{eqnarray*}

\noindent
Notice that for all $x\in S^2$, we have
\begin{equation}
\left( O,x,\left({x\over 2},{1\over 2}\right)\right)\in 
\Gamma^{\partial}_{O,x}\subset
\Gamma^{\mathcal U}_{O,x}\subset \Gamma^{\mathcal A}_{O,x} \label{eqn;a}
\end{equation}
\begin{equation}
d^{\Lambda,s}(O,x)\leq d^{s\left({x\over 2},{1\over 2}\right)}(O,x)=
{{1+f(x)}\over 2}\leq 1. 
\label{eqn;b}
\end{equation}

\indent
By the following two lemmas, we can use $\Gamma^{\partial}_{O,x}$ instead of
$\Gamma^{\mathcal{A}}_{O,x}$ in eq. (\ref{eqn;kdis}).

\begin{lem}\label{lem;adis} 
Suppose that $x\in S^2$. Then
$$ d^{\Lambda,s}(O,x)=\inf_{[X,\lambda,n]
\in \Gamma^{\mathcal U}_{O,x}}\sum^n_{i=1}d^{s(i)}(X_{i-1},X_i).$$ 
\end{lem}

\begin{proof}
Suppose that $[X,\lambda,n]\in \Gamma^{\mathcal A}_{O,x}\setminus 
\Gamma^{\mathcal U}_{O,x}$. If there is no $i$ such that 
$X_{i-1}$ and $X_i$ are antipodal on $S_{\lambda_i}$ with $\lambda_i=
\left({y\over 2},{1\over 2}\right)$ for some $y\in S^2$, then 
$d^{s(i)}(X_{i-1},X_i)=d_E(X_{i-1},X_i)$ for all $1\leq i\leq n$. Therefore  
$$\sum^n_{i=1}d^{s(i)}(X_{i-1},X_i)=\sum^n_{i=1}d_E(X_{i-1},X_i)
\geq d_E(O,x)=1.$$

\indent
If there exist two distinct integers $j,k$ such that $X_{j-1}$ and $X_j$
are antipodal on $S_{\lambda_j}$ with $\lambda_j=\left({y\over 2},{1\over 2}
\right)$, and $X_{k-1}$ and $X_k$ are antipodal on $S_{\lambda_k}$ with
$\lambda_k=\left({z\over 2},{1\over 2}\right)$ for some $y,z\in S^2$, 
then
\begin{eqnarray*}
\sum^n_{i=1}d^{s(i)}(X_{i-1},X_i)&\geq& 
d^{s(j)}(X_{j-1},X_j)+d^{s(k)}(X_{k-1},X_k) \\
&\geq&{{1+f(y)}\over 2}+{{1+f(z)}\over 2} \\
&\geq& 1. 
\end{eqnarray*}

\indent
Therefore by eq. (\ref{eqn;a}) and (\ref{eqn;b}), we have
\begin{eqnarray*}
d^{\Lambda,s}(O,x)=\inf_{[X,\lambda,n]\in \Gamma^{\mathcal U}_{O,x}}
\sum^n_{i=1}d^{s(i)}(X_{i-1},X_i). 
\end{eqnarray*}
\end{proof}

\begin{lem}\label{lem;f}
Suppose that $x\in S^2$. Then
\[ d^{\Lambda,s}(O,x)=\inf_{[X,\lambda,n]\in \Gamma^{\partial}_
{O,x}}\sum^n_{i=1}d^{s(i)}(X_{i-1},X_i). \]
\end{lem}

\begin{proof}
Suppose that $[X,\lambda,n]\in \Gamma^{\mathcal U}_{O,x}\setminus 
\Gamma^{\partial}_{O,x}$. There exists unique $j\neq 1$ such that 
$X_{j-1}$ and $X_j$ are antipodal on $S_{\lambda_j}$ with $\lambda_j=
\left({y\over 2},{1\over 2}\right)$ for some $y\in S^2$.
Let $X_{j-1}=z$ and $X_j=-z$ (see Figure \ref{fig;pf}).
Note that $z$ and $-z$ are antipodal on $S_{\lambda_j}$. Note also that $O$ 
and $y$ are antipodal on $S_{\lambda_j}$, too. 

\indent
If $d^{\Lambda,s}(O,z)\geq {1\over 2}$ then from eq. (\ref{eqn;a}) and  
(\ref{eqn;b}), we have  

\begin{eqnarray}
\sum^n_{i=1}d^{s(i)}(X_{i-1},X_i)&\geq& \sum^{j-1}_{i=1}
d^{s(i)}(X_{i-1},X_i)+d^{s\left({y\over 2},{1\over 2}\right)}(z,-z)\nonumber \\
&\geq& d^{\Lambda,s}(O,z)+d^{s\left({y\over 2},{1\over 2}\right)}(z,-z)
\nonumber \\
&\geq& {1\over 2}+{{1+f(y)}\over 2}\nonumber \\
&\geq& 1\nonumber\\ 
&\geq& {{1+f(x)}\over 2}\nonumber\\
&\geq& \inf_{[X,\lambda,n]\in \Gamma^{\partial}_{O,x}}\sum^n_{i=1}d^{s(i)}
(X_{i-1},X_i).
\label{eqn;p1} 
\end{eqnarray}

\indent
If $d^{\Lambda,s}(O,z)<{1\over 2}$ then from Lemma \ref{lem;half} and 
\ref{lem;haf}, we have 
$$d^{\Lambda,s}(O,z)=d_E(O,z)=d_E(y,-z)=d^{\Lambda,s}(y,-z).$$

\begin{figure}[ht]
\begin{center}
\psfrag{O}{$O$}
\psfrag{x}{$x$}
\psfrag{y}{$y$}
\psfrag{z}{$z$}
\psfrag{-z}{$-z$}
\includegraphics[width=2in, height=2in]{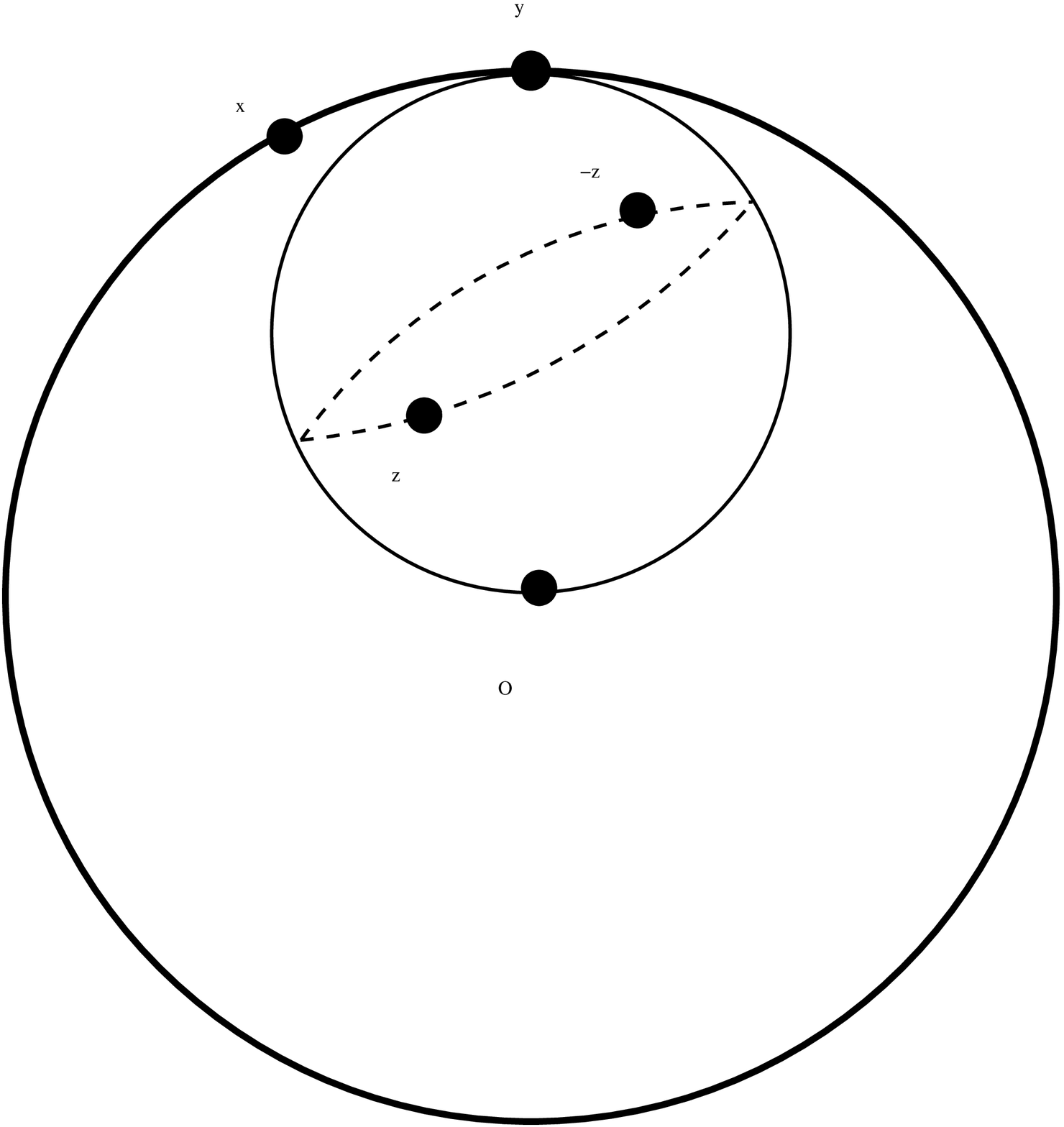}
\caption{$z\neq O$}
\label{fig;pf}
\end{center}
\end{figure}

\noindent
Therefore
\begin{eqnarray}
& &\sum^n_{i=1}d^{s(i)}(X_{i-1},X_i)\nonumber\\
&=&\sum^{j-1}_{i=1}d^{s(i)}(X_{i-1},X_i)
+d^{s\left({y\over 2},{1\over 2}\right)}(z,-z)+\sum^{n}_{i=j+1}d^{s(i)}
(X_{i-1},X_i)\nonumber\\
&\geq& d^{\Lambda,s}(O,z)+{{1+f(y)}\over 2}+d^{\Lambda,s}(-z,x)\nonumber \\
&=&d^{\Lambda,s}(y,-z)+{{1+f(y)}\over 2}+d^{\Lambda,s}(-z,x)\nonumber\\
&\geq&{{1+f(y)}\over 2}+d^{\Lambda,s}(y,x)\nonumber\\
&=& {{1+f(y)}\over 2}+\inf_{[X,\lambda,n]\in \Gamma_{y,x}}
\sum^n_{i=1}d^{s(i)}(X_{i-1},X_i).\label{eqn;p2}
\end{eqnarray}

\noindent
Suppose that $\left(Y_0,Y_1,\cdots,Y_m,\lambda'_1,\cdots,\lambda'_m\right)
\in\Gamma_{y,x}$ and let $s'(i)=s(\lambda'_i)$ for all $1\leq i\leq m$. 
If there exists $j$ such that 
$Y_{j-1}$ and $Y_j$ are antipodal on $S_{\lambda'_j}$ with $\lambda'_j=
\left({w\over 2},{1\over 2}\right)$ for some $w\in S^2$, then 
 from eq. (\ref{eqn;a}) and (\ref{eqn;b}) we have 
\begin{eqnarray*}
{{1+f(y)}\over 2}+\sum^m_{i=1}d^{s'(i)}(Y_{i-1},Y_i)\geq 
{{1+f(y)}\over 2}+d^{s'(j)}(Y_{j-1},Y_j)={{1+f(y)}\over 2}\qquad\\
+{{1+f(w)}\over 2}\geq 1\geq {{1+f(x)}\over 2}\geq\inf_{[X,\lambda,n]\in 
\Gamma^{\partial}_{O,x}}\sum^n_{i=1}d^{s(i)}(X_{i-1},X_i).
\end{eqnarray*} 
If there is no $i$ such that 
$Y_{i-1}$ and $Y_i$ are antipodal on $S_{\lambda'_i}$ with $\lambda'_i=
\left({w\over 2},{1\over 2}\right)$ for some $w\in S^2$, then from 
the definition of $\Gamma^{\partial}_{O,x}$ we have
$${{1+f(y)}\over 2}+\sum^m_{i=1}d^{s'(i)}(Y_{i-1},Y_i)\geq 
\inf_{[X,\lambda,n]\in \Gamma^{\partial}_{O,x}}\sum^n_{i=1}d^{s(i)}
(X_{i-1},X_i).$$

\noindent
Therefore in eq. (\ref{eqn;p2}) we have
$${{1+f(y)}\over 2}+\inf_{[X,\lambda,n]\in \Gamma_{y,x}}
\sum^n_{i=1}d^{s(i)}(X_{i-1},X_i)\geq 
\inf_{[X,\lambda,n]\in \Gamma^{\partial}_{O,x}}\sum^n_{i=1}d^{s(i)}
(X_{i-1},X_i).$$

\indent
Therefore from eq. (\ref{eqn;a}), (\ref{eqn;b}), (\ref{eqn;p1}) and 
Lemma \ref{lem;adis}, we have  
$$d^{\Lambda,s}(O,x)=\inf_{[X,\lambda,n]\in \Gamma^{\partial}_
{O,x}}\sum^n_{i=1}d^{s(i)}(X_{i-1},X_i).$$ 
\end{proof}

\indent
Due to Lemma \ref{lem;f} we can prove the following proposition which will be
used in the proof of Theorem \ref{thm;enough}. 

\begin{prop}\label{prop;nlf}
Suppose that $x\in S^2$. If $d^{\Lambda,s}(O,x)\neq 
{{1+f(x)}\over 2}$ then $x\in\mathcal{NM}^2_f$.
\end{prop}

\begin{proof}
Suppose that $x\in S^2$ and $d^{\Lambda,s}(O,x)\neq {{1+f(x)}\over 2}$.
From eq. (\ref{eqn;b}) we have 
$$d^{\Lambda,s}(O,x)<{{1+f(x)}\over 2}.$$
By Lemma \ref{lem;f}, there exists $y\in S^2$ such that
\begin{eqnarray*}
&&\left(O,y=X_1,X_2,\cdots,X_{n-1},x=X_n,\left({1\over 2}y,{1\over 2}
\right),\lambda_2,\cdots,\lambda_n\right)\in \Gamma^{\partial}_{O,x}\\ 
&&\qquad\qquad {{1+f(x)}\over 2}> d^{s\left({y\over 2},{1\over 2}\right)}(O,y)
+\sum^n_{i=2}d^{s(i)}(X_{i-1},X_i).
\end{eqnarray*}

\noindent
Since
$${{1+f(x)}\over 2}>d^{s\left({y\over 2},{1\over 2}\right)}(O,y)+\sum^n_{i=2}
d^{s(i)}(X_{i-1},X_i) 
\geq {{1+f(y)}\over 2}+d^{\Lambda,s}(y,x),$$
we have $1\geq f(x)-f(y)>2d^{\Lambda,s}(y,x)$.
Therefore from Lemma \ref{lem;half}. we have $d^{\Lambda,s}(y,x)=d_E(y,x)$.  
Thus $f(x)-f(y)>2d_E(x,y)$. Hence $x\in\mathcal{NM}^2_f$.
\end{proof}

\indent
We need the following fact in the proof of Theorem \ref{thm;enough}.

\begin{fact}\label{fct;path}
$\mb{R}^2\setminus \mathit C$ is path-connected if $\mathit C$ is a 
countable subset of $\mb{R}^2$. Therefore $S^2\setminus \mathit C$ 
is path-connected if $\mathit C$ is a countable subset of $S^2$.
\end{fact}

\noindent
Suppose that ${\mathit A}=\mb{R}^2\setminus \mathit{C}$ and $P,Q\in\mathit A$.
There are uncountable straight lines in $\mathit{A}$ which contain $P$ and 
there are uncountable straight lines in $\mathit{A}$ which contain $Q$.
We can choose non-parallel two lines to find a path from $P$ to $Q$.  

\indent 
Now we can prove Theorem \ref{thm;enough}.

\begin{pfenough}
Suppose that $f:S^2 \rightarrow [0,1]$ is a countable-to-one function.
To get a contradiction suppose that $\mathcal{NM}^2_f$ is countable.
Notice that $\bigcup_{t\in \mb{Q}}f^{-1}(t)$ is countable, where $\mb{Q}$
is the set of rational numbers. Therefore the following set ${\mathit A}$ is 
path-connected by Fact \ref{fct;path}.
$${\mathit A}=S^2\setminus \left(\left(\bigcup_{t\in\mb{Q}}f^{-1}(t)
\right)\cup \mathcal{NM}^2_f\right)$$
Define a function $g:({\mathit A},d_E)\to [0,1]$ by $g(x)=
d^{\Lambda,s}(O,x)$. Since $\Lambda$ is piecewise dense, we have 
$$ |g(x)-g(y)|=|d^{\Lambda,s}(O,x)-d^{\Lambda,s}(O,y)|\leq d^{\Lambda,s}
(x,y)\leq d_E(x,y).$$
Therefore $g$ is a continuous function and hence $g({\mathit A})$ is 
path-connected.

\indent
From Proposition \ref{prop;nlf}, we have $g(x)={{1+f(x)}\over 2}$
for all $x\in {\mathit A}$. Therefore $g(x)$ is an irrational number for all
$x\in{\mathit A}$. Notice that $g^{-1}(t)$ is countable for all $t$ 
and $\mathit A$ is a uncountable set. Therefore $g({\mathit A})$ is a 
uncountable subset of irrational numbers. Therefore $g({\mathit A})$ is not 
path-connected. This is a contradiction. Therefore $\mathcal{NM}^2_f$ is 
uncountable. 
\end{pfenough}

\end{document}